\documentclass[12pt]{article}
\usepackage{latexsym, amssymb, amsmath, amscd, amsfonts, epsfig, graphicx, colordvi,verbatim,ifpdf,extarrows}
\usepackage{amsfonts, amsmath, amssymb}
\usepackage{amssymb,amsfonts,amsmath,latexsym,epsfig,cite, psfrag,eepic,color}
\usepackage{amscd,graphics}
\usepackage{latexsym, amssymb,  amsmath,amscd, amsfonts, epsfig, graphicx, colordvi,amsthm}

\usepackage{graphicx}
\usepackage{epstopdf}
\usepackage{color}
\usepackage{ifpdf}
\usepackage{fancybox}
\usepackage[font=small,labelfont=bf,labelsep=none]{caption}
\usepackage{float}

\allowdisplaybreaks

\usepackage[latin1]{inputenc}

\newtheorem{theorem}{Theorem}[section]

\theoremstyle{remark}

\theoremstyle{definition}

\newtheorem{core}[theorem]{Corollary}

\def\qed{\nopagebreak\hfill{\rule{4pt}{7pt}}
\medbreak}

\numberwithin{equation}{section}

\def\qed{\nopagebreak\hfill{\rule{4pt}{7pt}}
\medbreak}

\newlength{\boxedparwidth}
  {\begin{center} \begin{tabular}{|@{\hspace{.315in}}c@{\hspace{.15in}}|}
                  \hline \\ \begin{minipage}[t]{\boxedparwidth}
                  \setlength{\parindent}{.25in}}%
  {\end{minipage} \\ \\ \hline \end{tabular} \end{center}}

\begin{document}

\begin{center}
{\Large \bf Three types of the minimal excludant size of an overpartition}
\end{center}

\begin{center}
 {Thomas Y. He}$^{1}$, {C.S. Huang}$^{2}$,
  {H.X. Li}$^{3}$ and {X. Zhang}$^{4}$ \vskip 2mm

$^{1,2,3,4}$ School of Mathematical Sciences, Sichuan Normal University, Chengdu 610066, P.R. China

   \vskip 2mm

  $^1$heyao@sicnu.edu.cn,  $^2$huangchushu@stu.sicnu.edu.cn,  $^3$lihaixia@stu.sicnu.edu.cn, $^4$zhangxi@stu.sicnu.edu.cn
\end{center}

\vskip 6mm   {\noindent \bf Abstract.} Recently, Andrews and Newman studied the minimal excludant of a partition, which is defined as the smallest positive integer that is not a part of a partition. In this article, we consider the minimal excludant size of an overpartition, which is an extension of  the minimal excludant of a partition. We will investigate three types of the minimal excludant size of an overpartition.

\section{Introduction}

A partition $\pi$ of a positive integer $n$ is a finite non-increasing sequence of positive integers $\pi=(\pi_1,\pi_2,\ldots,\pi_\ell)$ such that $\pi_1+\pi_2+\cdots+\pi_\ell=n$. The $\pi_i$ are called the parts of $\pi$. The weight of $\pi$ is the sum of parts, denoted $|\pi|$.

The minimal excludant of a partition was introduced by Grabner and Knopfmacher \cite{Grabner-Knopfmacher-2006} under the name ``smallest gap". Recently, Andrews and Newman \cite{Andrews-Newman-2019} undertook a combinatorial study of the minimal excludant of a partition. The minimal excludant of a partition $\pi$ is the smallest positive integer that is not a part of $\pi$, denoted  $mex(\pi)$.

There are two generalization of the minimal excludant of a partition. In \cite{Andrews-Newman-2020}, Andrews and Newman defined $mex_{A,a}(\pi)$ to be the smallest positive integer congruent to $a$ modulo $A$ that is not a part of $\pi$. In \cite{Ballantine-Merca-2020}, Ballantine and Merca introduced the $r$-gaps of a partition $\pi$, denoted $g_r(\pi)$, which is the smallest positive integer that does not appear at least $r$ times in $\pi$. Clearly, \[mex_{1,1}(\pi)=mex(\pi)\text{ and }g_1(\pi)=mex(\pi).\]

In this article, we extended the minimal excludant of partitions to overpartitions involving $mex_{A,a}(\pi)$ and $g_r(\pi)$. An overpartition, introduced by Corteel and Lovejoy \cite{Corteel-Lovejoy-2004},  is a partition such that the first occurrence of a part can be overlined.
For example, there are fourteen overpartitions of $4$:
\[(4),(\bar{4}),(3,1),(\bar{3},1),(3,\bar{1}),(\bar{3},\bar{1}),(2,2),(\bar{2},2),
\]
\[(2,1,1),(\bar{2},1,1),(2,\bar{1},1),(\bar{2},\bar{1},1),(1,1,1,1),(\bar{1},1,1,1).\]

Throughout this article, we adopt the following convention.  We use $\overline{\mathcal{P}}$ to denote the set of all overpartitions and use $\overline{\mathcal{P}}(n)$ to denote the set of all overpartitions of $n$.  Let $\overline{p}(n)$ be the number of overpartitions of $n$. Clearly, we have $\overline{p}(4)=14$. The generating function of $\overline{p}(n)$ is
\begin{equation}\label{gen-over-p}
\sum_{n=0}^\infty\overline{p}(n)q^n=\frac{(-q;q)_\infty}{(q;q)_\infty}.
\end{equation}
Here and in the sequel, we assume that $|q|<1$ and use the standard notation \cite{Andrews-1976}:
\[(a;q)_\infty=\prod_{i=0}^{\infty}(1-aq^i)\quad\text{and}\quad(a;q)_n=\frac{(a;q)_\infty}{(aq^n;q)_\infty}.\]

Let $\pi$ be an overpartition. For a part $\pi_i$ of $\pi$, we say that $\pi_i$ is of size $t$ if $\pi_i=t$ or $\overline{t}$. We use $f_t(\pi)$ (resp. $f_{\bar{t}}(\pi)$) to denote the number of parts equal to $t$ (resp. $\bar{t}$) in $\pi$. It follows from the definition of overpartition that $f_{\bar{t}}(\pi)=0$ or $1$.

Aricheta and Donato \cite{Aricheta-Donato-2023}, Yang and Zhou \cite{Yang-Zhou-2023,Yang-Zhou-2024} and Dhar, Mukhopadhyay and Sarma \cite{Dhar-Mukhopadhyay-Sarma-2023,Dhar-Mukhopadhyay-Sarma-2024} have extended the minimal excludant of partitions to overpartitions. In this article, we introduce three types of the minimal excludant size of an overpartition, which are stated as follows. In the remaining of this article, we assume that $A\geq a\geq 1$.

\begin{itemize}
\item[(1)] For $r\geq1$, let $mes_{r,A,a}(\pi)$ be the smallest positive integer $\equiv a\pmod A$ such that there are no $r$ parts of size $mes_{r,A,a}(\pi)$ in $\pi$.

\item[(2)] For $r\geq2$, let $\overline{mes}_{r,A,a}(\pi)$ be the smallest positive integer $\equiv a\pmod A$ such that there is no overlined part of size $\overline{mes}_{r,A,a}(\pi)$ in $\pi$ or there are no $r-1$ non-overlined parts of size $\overline{mes}_{r,A,a}(\pi)$ in $\pi$.

\item[(3)] For  $r\ge$ 2, let $\widetilde{mes}_{r,A,a}(\pi)$ be the smallest positive integer $\equiv a\pmod A$ such that there is no overlined part of size $\widetilde{mes}_{r,A,a}(\pi)$ in $\pi$ and there are no $r-1$ non-overlined parts of size $\widetilde{mes}_{r,A,a}(\pi)$ in $\pi$.

\end{itemize}

The minimal excludant size $mes_{1,1,1}(\pi)$ of an overpartition $\pi$ is the $\overline{mex}(\pi)$ in \cite{Dhar-Mukhopadhyay-Sarma-2023,Dhar-Mukhopadhyay-Sarma-2024}. For example, the $mes_{r,A,a}(\pi)$, $\overline{mes}_{r,A,a}(\pi)$ and $\widetilde{mes}_{r,A,a}(\pi)$ of the overpartitions $\pi$ of $4$ are given in the following table.

\begin{center}
\begin{tabular}{|c|c|c|c|c|c|c|c|c|}
  \hline
  $\pi$&$(4)$&$(\bar{4})$&$(3,1)$&$(\bar{3},1)$&$(3,\bar{1})$&$(\bar{3},\bar{1})$&$(2,2)$&$(\bar{2},2)$\\
  \hline
  ${mes}_{2,2,1}(\pi)$&$1$&$1$&$1$&$1$&$1$&$1$&$1$&$1$\\
  \hline
  $\overline{mes}_{2,2,1}(\pi)$&$1$&$1$&$1$&$1$&$1$&$1$&$1$&$1$\\
  \hline
  $\widetilde{mes}_{2,2,1}(\pi)$&$1$&$1$&$5$&$5$&$5$&$5$&$1$&$1$\\
  \hline
\end{tabular}
\end{center}

\begin{center}
\begin{tabular}{|c|c|c|c|c|c|c|}
  \hline
  $\pi$&$(2,1,1)$&$(\bar{2},1,1)$&$(2,\bar{1},1)$&$(\bar{2},\bar{1},1)$&$(1,1,1,1)$&$(\bar{1},1,1,1)$\\
  \hline
  ${mes}_{2,2,1}(\pi)$&$3$&$3$&$3$&$3$&$3$&$3$\\
  \hline
   $\overline{mes}_{2,2,1}(\pi)$&$1$&$1$&$3$&$3$&$1$&$3$\\
  \hline
  $\widetilde{mes}_{2,2,1}(\pi)$&$3$&$3$&$3$&$3$&$3$&$3$\\
  \hline
\end{tabular}
\end{center}

This article is organized as follows. In Section 2, we will list the results of this article in Section 2. We then give the proofs of the results related to $mes_{r,A,a}(\pi)$, $\overline{mes}_{r,A,a}(\pi)$ and $\widetilde{mes}_{r,A,a}(\pi)$ in Section 3, Section 4 and Section 5 respectively.

\section{Main results of this article}

In this section, we will list the results related to $mes_{r,A,a}(\pi)$, $\overline{mes}_{r,A,a}(\pi)$ and $\widetilde{mes}_{r,A,a}(\pi)$ in Section 2.1, Section 2.2 and Section 2.3 respectively.

\subsection{Results related to $mes_{r,A,a}(\pi)$}

We can obtain the generating function for all overpartitions by considering $mes_{r,A,a}(\pi)$.
\begin{theorem}\label{gen-mes}
For $r\geq 1$,
\begin{equation*}\label{gen-mes-eqn}
\sum_{\pi\in\mathcal{\overline{P}}}z^{mes_{r,A,a}(\pi)}q^{|\pi|}
=\frac{(-q;q)_{\infty}}{(q;q)_{\infty}}\sum_{k=0}^{\infty}z^{kA+a}\bigg[\frac{2^kq^{r[A\binom{k}{2}+ka]}}{(-q^a;q^A)_k}-\frac{2^{k+1}q^{r[A\binom{k+1}{2}+(k+1)a]}}{(-q^a;q^A)_{k+1}}\bigg].
\end{equation*}
\end{theorem}

For $r\geq 1$ and $n\geq 0$, we define \[\sigma mes_{r,A,a}(n)=\sum_{\pi\in \overline{\mathcal{P}}(n)}mes_{r,A,a}(\pi).\]
 If $r=A=a=1$, then $\sigma mes_{1,1,1}(n)$ is the $\sigma \overline{mex}(n)$ in \cite{Dhar-Mukhopadhyay-Sarma-2024}. The generating function of $\sigma mes_{r,A,a}(n)$ is given as follows.
\begin{theorem}\label{gen-sigma-mes}For $r\geq 1$,
\begin{equation*}\label{gen-sigma-mes-eqn}
\sum_{n=0}^\infty\sigma mes_{r,A,a}(n)q^n
=\frac{(-q;q)_{\infty}}{(q;q)_{\infty}}\bigg[a+A\sum_{k=1}^{\infty}\frac{2^kq^{r[A\binom{k}{2}+ka]}}{(-q^a;q^A)_k}\bigg].
\end{equation*}
\end{theorem}

For $r\geq 1$, let $Nmes_{r,A,a}$ be set of overpartitions $\pi$ such that all parts $\equiv a\pmod A$ of size less than  $mes_{r,A,a}(\pi)$ in $\pi$ are non-overlined. If $r=A=a=1$, then the number of overpartitions of $n$ in $Nmes_{1,1,1}$ is the $\overline{m}(n)$ in \cite{Dhar-Mukhopadhyay-Sarma-2023,Dhar-Mukhopadhyay-Sarma-2024}.

\begin{theorem}\label{gen-Nmes} For $r\geq 1$,
\begin{equation}\label{gen-Nmes-eqn}
\begin{split}
\sum_{\pi\in Nmes_{r,A,a}}z^{mes_{r,A,a}(\pi)}q^{|\pi|}=\frac{(-q;q)_{\infty}}{(q;q)_{\infty}}\sum_{k=0}^{\infty}z^{kA+a}\bigg[\frac{q^{r[A\binom{k}{2}+ka]}}{(-q^a;q^A)_k}-\frac{2
q^{r[A\binom{k+1}{2}+(k+1)a]}}{(-q^a;q^A)_{k+1}}\bigg].
\end{split}
\end{equation}
\end{theorem}

Setting $z=1$ in \eqref{gen-Nmes-eqn}, we can get
\begin{equation}\label{gen-n-cor}
\sum_{\pi\in Nmes_{r,A,a}}q^{|\pi|}=
\frac{(-q;q)_{\infty}}{(q;q)_{\infty}}-\frac{(-q;q)_{\infty}}{(q;q)_{\infty}}\sum_{k=1}^{\infty}\frac{q^{r[A\binom{k}{2}+ka]}}{(-q^a;q^A)_k}.
\end{equation}
Combining with \eqref{gen-over-p}, we can obtain the following corollary.
\begin{core}
The generating function for the overpartitions $\pi$ such that there exist overlined parts $\equiv a\pmod A$ of size less than $mes_{r,A,a}(\pi)$ in $\pi$ is
\[\frac{(-q;q)_{\infty}}{(q;q)_{\infty}}\sum_{k=1}^{\infty}\frac{q^{r[A\binom{k}{2}+ka]}}{(-q^a;q^A)_k}.\]
\end{core}

Letting $r=A=a=1$ in \eqref{gen-n-cor}, we get
\begin{align*}
\sum_{\pi\in Nmes_{1,1,1}}q^{|\pi|}&=
\frac{(-q;q)_{\infty}}{(q;q)_{\infty}}-\frac{(-q;q)_{\infty}}{(q;q)_{\infty}}\sum_{k=1}^{\infty}\frac{q^{\binom{k+1}{2}}}{(-q;q)_k}\\
&=\frac{(-q;q)_{\infty}}{(q;q)_{\infty}}-\frac{(-q;q)_{\infty}}{(q;q)_{\infty}}\sum_{k=1}^{\infty}\frac{q^{\binom{k}{2}}(q^k+1-1)}{(-q;q)_k}\\
&=\frac{(-q;q)_{\infty}}{(q;q)_{\infty}}-\frac{(-q;q)_{\infty}}{(q;q)_{\infty}}\sum_{k=1}^{\infty}\frac{q^{\binom{k}{2}}}{(-q;q)_{k-1}}+\frac{(-q;q)_{\infty}}{(q;q)_{\infty}}\sum_{k=1}^{\infty}\frac{q^{\binom{k}{2}}}{(-q;q)_{k}}\\
&=\frac{(-q;q)_{\infty}}{(q;q)_{\infty}}\sum_{k=0}^{\infty}\frac{q^{\binom{k}{2}}}{(-q;q)_{k}}-\frac{(-q;q)_{\infty}}{(q;q)_{\infty}}\sum_{k=0}^{\infty}\frac{q^{\binom{k+1}{2}}}{(-q;q)_{k}},
\end{align*}
which leads to Theorem 30 in \cite{Dhar-Mukhopadhyay-Sarma-2023}.

For $r=1$, let $Omes_{1,A,a}$ be the set of overpartitions $\pi$ such that all parts $\equiv a\pmod A$ of size less than $mes_{1,A,a}(\pi)$ in $\pi$ are overlined.
 \begin{theorem}\label{gen-Omes} The generating function for the overpartitions in $Omes_{1,A,a}$ is
\begin{equation}\label{gen-Omes-eqn}
\begin{split}
&\sum_{\pi\in Omes_{1,A,a}}z^{mes_{1,A,a}(\pi)}q^{|\pi|}\\
&=\frac{(-q;q)_{\infty}}{(q;q)_{\infty}}\sum_{k=0}^{\infty}z^{kA+a}\bigg[q^{A\binom{k}{2}+ka} \frac{(q^a;q^A)_k}{(-q^a;q^A)_k}-2q^{A\binom{k+1}{2}+(k+1)a} \frac{(q^a;q^A)_k}{(-q^a;q^A)_{k+1}}\bigg].
\end{split}
\end{equation}
\end{theorem}

Setting $z=1$ in \eqref{gen-Omes-eqn}, we can get
\begin{align*}
\sum_{\pi\in Omes_{1,A,a}}q^{|\pi|}&=\frac{(-q;q)_{\infty}}{(q;q)_{\infty}}\bigg[\sum_{k=0}^{\infty}q^{A\binom{k}{2}+ka} \frac{(q^a;q^A)_k}{(-q^a;q^A)_k}-
2\sum_{k=1}^{\infty}q^{A\binom{k}{2}+ka} \frac{(q^a;q^A)_{k-1}}{(-q^a;q^A)_k}\bigg]\nonumber\\
&=\frac{(-q;q)_{\infty}}{(q;q)_{\infty}}-\frac{(-q;q)_{\infty}}{(q;q)_{\infty}}\sum_{k=1}^{\infty}q^{A\binom{k}{2}+ka} \frac{(q^a;q^A)_{k-1}}{(-q^a;q^A)_k}[2-(1-q^{(k-1)A+a})]\\
&=\frac{(-q;q)_{\infty}}{(q;q)_{\infty}}-\frac{(-q;q)_{\infty}}{(q;q)_{\infty}}\sum_{k=1}^{\infty}q^{A\binom{k}{2}+ka} \frac{(q^a;q^A)_{k-1}}{(-q^a;q^A)_{k-1}}.
\end{align*}

We can obtain the following corollary.
\begin{core}
The generating function for the overpartitions $\pi$ such that there exist non-overlined $parts\equiv a\pmod A$ of size less than $mes_{1,A,a}(\pi)$ in $\pi$ is
\[\frac{(-q;q)_{\infty}}{(q;q)_{\infty}}\sum_{k=1}^{\infty}q^{A\binom{k}{2}+ka} \frac{(q^a;q^A)_{k-1}}{(-q^a;q^A)_{k-1}}.\]
\end{core}

 We will give the proofs of Theorems \ref{gen-mes}, \ref{gen-sigma-mes}, \ref{gen-Nmes} and \ref{gen-Omes} in Section 3.

\subsection{Results related to $\overline{mes}_{r,A,a}(\pi)$}

We can obtain the generating function for all overpartitions in terms of $\overline{mes}_{r,A,a}(\pi)$.
\begin{theorem}\label{gen-overmes}
For $r\geq 2$,
\begin{equation*}\label{gen-overmes-eqn}
\sum_{\pi\in\mathcal{\overline{P}}}z^{\overline{mes}_{r,A,a}(\pi)}q^{|\pi|}
=\frac{(-q;q)_{\infty}}{(q;q)_{\infty}}\sum_{k=0}^{\infty}z^{kA+a}\bigg[\frac{q^{r[A\binom{k}{2}+ka]}}{(-q^a;q^A)_k}-\frac{q^{r[A\binom{k+1}{2}+(k+1)a]}}{(-q^a;q^A)_{k+1}}\bigg].
\end{equation*}
\end{theorem}

For $r\geq 2$ and $n\geq 0$, we define
\[
\sigma\overline{mes}_{r,A,a}(n)=\sum_{\pi\in\overline{\mathcal{P}}(n)}\overline{mes}_{r,A,a}(\pi).
\]
The generating function of $\sigma\overline{mes}_{r,A,a}(n)$ is given as follows.
\begin{theorem}\label{gen-sigma-overmes}For $r\geq 2$,
\begin{equation*}\label{gen-sigma-overmes-eqn}
\sum_{n=0}^\infty\sigma\overline{mes}_{r,A,a(n)}q^{n}
=\frac{{(-q;q)_{\infty}}}{{(q;q)_{\infty}}}\bigg[a+A\sum_{k=1}^{\infty}\frac{q^{r[A\binom{k}{2}+ka]}}{(-q^{a};q^{A})_{k}}\bigg].
\end{equation*}
\end{theorem}

For $r\geq 2$, let $N\overline{mes}_{r,A,a}$ (resp. $O\overline{mes}_{r,A,a}$) be the set of overpartition $\pi$ such that there are
no $r-1$ (resp. there are at least $r-1$) non-overlined parts of size $\overline{mes}_{r,A,a}(\pi)$ in $\pi$. If $r=2$ and $A=a=1$, then the number of overpartitions of $n$ in $O\overline{mes}_{2,1,1}$ is the $\overline{m}^{ord}(n)$ in \cite{Dhar-Mukhopadhyay-Sarma-2024}.

\begin{theorem}\label{gen-overNmes} For $r\geq 2$,
\begin{equation}\label{gen-overNmes-eqn}
\sum_{\pi\in N\overline{mes}_{r,A,a}}z^{\overline{mes}_{r,A,a}(\pi)}q^{|\pi|}=\frac{{(-q;q)_{\infty}}}{{(q;q)_{\infty}}}\sum_{k=0}^{\infty}z^{kA+a}\bigg[\frac{q^{r[A\binom{k}{2}+ka]}}{(-q^{a};q^{A})_{k}}-\frac{q^{r[A\binom{k}{2}+ka]+(r-1)(kA+a)}}{(-q^{a};q^{A})_{k}}\bigg].
\end{equation}
\end{theorem}
Setting $z=1$ in \eqref{gen-overNmes-eqn}, we can get
\begin{align*}
&\sum_{\pi\in N\overline{mes}_{r,A,a}}q^{|\pi|}\\
&=\frac{{(-q;q)_{\infty}}}{{(q;q)_{\infty}}}\sum_{k=0}^{\infty}\bigg[\frac{q^{r[A\binom{k}{2}+ka]}}{(-q^{a};q^{A})_{k}}-\frac{q^{r[A\binom{k}{2}+ka]+(r-1)(kA+a)}}{(-q^{a};q^{A})_{k}}\frac{1+q^{kA+a}}{1+q^{kA+a}}\bigg]\\
&=\frac{{(-q;q)_{\infty}}}{{(q;q)_{\infty}}}\bigg[\sum_{k=0}^{\infty}\frac{q^{r[A\binom{k}{2}+ka]}}{(-q^{a};q^{A})_{k}}-\sum_{k=0}^{\infty}\frac{q^{r[A\binom{k}{2}+ka]+(r-1)(kA+a)}}{(-q^{a};q^{A})_{k+1}}-\sum_{k=0}^{\infty}\frac{q^{r[A\binom{k+1}{2}+(k+1)a]}}{(-q^{a};q^{A})_{k+1}}\bigg]\\
&=\frac{{(-q;q)_{\infty}}}{{(q;q)_{\infty}}}-\frac{{(-q;q)_{\infty}}}{{(q;q)_{\infty}}}\sum_{k=0}^{\infty}\frac{q^{r[A\binom{k}{2}+ka]+(r-1)(kA+a)}}{(-q^{a};q^{A})_{k+1}}.
\end{align*}

As a corollary, we can obtain the generating function for the overpartitions in $O\overline{mes}_{r,A,a}$.
\begin{core}\label{gen-O}For $r\geq 2$,
\begin{equation*}
\sum_{\pi\in O\overline{mes}_{r,A,a}}q^{|\pi|}
=\frac{{(-q;q)_{\infty}}}{{(q;q)_{\infty}}}\sum_{k=0}^{\infty}\frac{q^{r[A\binom{k}{2}+ka]+(r-1)(kA+a)}}{(-q^{a};q^{A})_{k+1}}.
\end{equation*}
\end{core}
Setting $r=2$ and $A=a=1$ in Corollary \ref{gen-O}, we get Theorem 2.3 in \cite{Dhar-Mukhopadhyay-Sarma-2024}.

The proofs of Theorems \ref{gen-overmes}, \ref{gen-sigma-overmes} and \ref{gen-overNmes} will be given in Section 4.

\subsection{Results related to $\widetilde{mes}_{r,A,a}(\pi)$}

We can obtain the generating function for all overpartitions in light of $\widetilde{mes}_{r,A,a}(\pi)$.
\begin{theorem}\label{gen-tildemes}
For $r\geq 2$,
\begin{equation*}\label{gen-tildemes-eqn}
\begin{split}
&\sum_{\pi\in\overline{P}}z^{\widetilde{mes}_{r,A,a}(\pi)} q^{|\pi|}\\
&=\frac{(-q;q)_{\infty }}{(q;q)_{\infty}}\sum_{k=0}^{\infty }z^{kA+a}\bigg[\frac{q^{A[{k\choose 2}+ka ]}}{(-q^{a};q^{A})_{k}}(-q^{(r-2)a};q^{(r-2)A})_{k}\\
&\qquad\qquad\qquad\qquad\qquad-\frac{q^{A[{k+1\choose 2}+(k+1)a ]}}{(-q^{a};q^{A})_{k+1}}(-q^{(r-2)a};q^{(r-2)A})_{k+1}\bigg].
\end{split}
\end{equation*}
\end{theorem}

For $r\geq 2$ and $n\geq 0$, we define
\[\sigma\widetilde{mes}_{r,A,a}(n)=\sum_{\pi\in\overline{\mathcal{P}}(n)}\widetilde{mes}_{r,A,a}(\pi).\]
The generating function of $\sigma\widetilde{mes}_{r,A,a}(n)$ is given as follows.
\begin{theorem}\label{gen-sigma-tildemes}For $r\geq 2$,
\begin{equation*}\label{gen-sigma-tildemes-eqn}
\sum_{n=0}^\infty\sigma\widetilde{mes}_{r,A,a}(n)q^{n}
=\frac{(-q;q)_{\infty }}{(q;q)_{\infty }}\bigg[a+A\sum_{k=1}^{\infty }\frac{q^{A[{k\choose 2}+ka]}}{(-q^{a};q^{A})_{k}}(-q^{(r-2)a};q^{(r-2)A})_{k}\bigg].
\end{equation*}
\end{theorem}
For $r\geq 2$, let $N\widetilde{mes}_{r,A,a}$ (resp. $O\widetilde{mes}_{r,A,a}$) be the set of overpartions $\pi$ such that all parts $\equiv a\pmod A$ of size less than $\widetilde{mes}_{r,A,a}(\pi)$ in $\pi$ are non-overlined (resp.  overlined).

\begin{theorem}\label{gen-tildeNmes} For $r\geq 2$,
\begin{equation}\label{gen-tildeNmes-eqn}
\sum_{\pi \in N\widetilde{mes}_{r,A,a}}z^{\widetilde{mes}_{r,A,a}(\pi)}q^{|\pi|}=\frac{(-q;q)_{\infty }}{(q;q)_{\infty}}\sum_{k=0}^{\infty }z^{kA+a}\bigg[\frac{q^{(r-1)[A{k\choose 2}+ka]}}{(-q^{a};q^{A})_{k+1 }}-\frac{q^{(r-1)[A{k+1\choose 2}+(k+1)a]}}{(-q^{a};q^{A})_{k+1}}\bigg].
\end{equation}
\end{theorem}

\begin{theorem}\label{gen-tildeOmes} For $r\geq 2$,
\begin{equation}\label{gen-tildeOmes-eqn}
\begin{split}
&\sum_{\pi \in O\widetilde{mes}_{r,A,a}}z^{\widetilde{mes}_{r,A,a}(\pi)}q^{|\pi|}\\
&=\frac{(-q;q)_{\infty }}{(q;q)_{\infty}}\sum_{k=0}^{\infty }z^{kA+a}\bigg[q^{A{k\choose 2}+ka}\frac{(q^{a};q^{A})_{k}}{(-q^{a};q^{A})_{k+1 }}-q^{A{k\choose 2}+ka+(r-1)(kA+a)}\frac{(q^{a};q^{A})_{k}}{(-q^{a};q^{A})_{k+1}}\bigg].
\end{split}
\end{equation}
\end{theorem}

Setting $z=1$ in \eqref{gen-tildeNmes-eqn}, we can get
\begin{align*}
	&\sum_{\pi \in N\widetilde{mes}_{r,A,a} }q^{|\pi|}\\
	&=\frac{(-q;q)_{\infty }}{(q;q)_{\infty }} \sum_{k=0}^{\infty } \bigg[\frac{q^{(r-1)[A{k\choose 2}+ka]}}{(-q^{a};q^{A})_{k+1}}
	(1+q^{kA+a}-q^{kA+a}) -\frac{q^{(r-1)[A{k+1\choose 2}+(k+1)a]}}{(-q^{a};q^{A})_{k+1}}\bigg]\\
	&= \frac{(-q;q)_{\infty }}{(q;q)_{\infty }}\bigg[\sum_{k=0}^{\infty } \frac{q^{(r-1)[A{k\choose 2}+ka]}}{(-q^{a};q^{A})_{k}}-\sum_{k=0}^{\infty }\frac{q^{(r-1)[A{k\choose 2}+ka]+kA+a}}{(-q^{a};q^{A})_{k+1}}
	-\sum_{k=1}^{\infty }  \frac{q^{(r-1)[A{k\choose 2}+ka]}}{(-q^{a};q^{A})_{k}}\bigg]\\
	&=\frac{(-q;q)_{\infty }}{(q;q)_{\infty }}-\frac{(-q;q)_{\infty }}{(q;q)_{\infty }}\sum_{k=0}^{\infty}\frac{q^{(r-1)[A{k\choose 2}+ka]+kA+a}}{(-q^{a};q^{A})_{k+1}}.
\end{align*}
We can obtain the following corollary.
\begin{core}
The generating function for the overpartitions $\pi$ such that there exist overlined parts $\equiv a \pmod A$ of size less than $\widetilde{mes}_{r,A,a}(\pi)$ in $\pi$ is
\[\frac{(-q;q)_{\infty }}{(q;q)_{\infty }}\sum_{k=0}^{\infty}\frac{q^{(r-1)[A{k\choose 2}+ka]+kA+a}}{(-q^{a};q^{A})_{k+1}}.\]
\end{core}

Setting $z=1$ in \eqref{gen-tildeOmes-eqn}, we can get

	\begin{align*}
	&\sum_{\pi \in O\widetilde{mes}_{r,A,a} }q^{|\pi|}\\
&=\frac{(-q;q)_{\infty }}{(q;q)_{\infty }} \bigg[\sum_{k=0}^{\infty } q^{A{k\choose 2}+ka}\frac{(q^{a};q^{A})_{k}}{(-q^{a};q^{A})_{k+1}}
	(1+q^{kA+a}-q^{kA+a})\\
	&\qquad\qquad\qquad-\sum_{k=0}^{\infty}q^{A{{k+1}\choose 2}+(k+1)a+(r-2)(kA+a)}\frac{(q^{a};q^{A})_{k}}{(-q^{a};q^{A})_{k+1}}\bigg]\\
	&=\frac{(-q;q)_{\infty }}{(q;q)_{\infty }} \bigg[\sum_{k=0}^{\infty } q^{A{k\choose 2}+ka}\frac{(q^{a};q^{A})_{k}}{(-q^{a};q^{A})_{k}}-\sum_{k=0}^{\infty}q^{A{k+1\choose 2}+(k+1)a}\frac{(q^{a};q^{A})_{k}}{(-q^{a};q^{A})_{k+1}}\\
	&\qquad\qquad\qquad-\sum_{k=1}^{\infty}q^{A{k\choose 2}+ka+(r-2)((k-1)A+a)}\frac{(q^{a};q^{A})_{k-1}}{(-q^{a};q^{A})_{k}}\bigg]\\
	&=\frac{(-q;q)_{\infty }}{(q;q)_{\infty }}
	\bigg[1+\sum_{k=1}^{\infty } q^{A{k\choose 2}+ka}\frac{(q^{a};q^{A})_{k}}{(-q^{a};q^{A})_{k}}-\sum_{k=1}^{\infty } q^{A{k\choose 2}+ka}\frac{(q^{a};q^{A})_{k-1}}{(-q^{a};q^{A})_{k}}\\
	&\qquad\qquad\qquad-\sum_{k=1}^{\infty}q^{A{k\choose 2}+ka+(r-2)((k-1)A+a)}\frac{(q^{a};q^{A})_{k-1}}{(-q^{a};q^{A})_{k}}\bigg]\\
	&=\frac{(-q;q)_{\infty }}{(q;q)_{\infty }}\bigg[1+\sum_{k=1}^{\infty } q^{A{k\choose 2}+ka}\frac{(q^{a};q^{A})_{k-1}}{(-q^{a};q^{A})_{k}}(1-q^{(k-1)A+a}-1-q^{(r-2)((k-1)A+a)})\bigg]\\
	&=\frac{(-q;q)_{\infty }}{(q;q)_{\infty }}-\frac{(-q;q)_{\infty }}{(q;q)_{\infty }}\sum_{k=1}^{\infty } q^{A{k\choose 2}+ka}\frac{(q^{a};q^{A})_{k-1}}{(-q^{a};q^{A})_{k}}(q^{(k-1)A+a}+q^{(r-2)((k-1)A+a)}),
	\end{align*}
which leads to the following corollary.
\begin{core}
The generating function for the overpartitions $\pi$ such that there exist non-overlined parts $\equiv a \pmod A$ of size less than $\widetilde{mes}_{r,A,a}(\pi)$ in $\pi$ is
\[\frac{(-q;q)_{\infty }}{(q;q)_{\infty }}\sum_{k=1}^{\infty } q^{A{k\choose 2}+ka}\frac{(q^{a};q^{A})_{k-1}}{(-q^{a};q^{A})_{k}}(q^{(k-1)A+a}+q^{(r-2)((k-1)A+a)}).\]
\end{core}

 We will present the proofs of Theorems \ref{gen-tildemes}, \ref{gen-sigma-tildemes}, \ref{gen-tildeNmes} and \ref{gen-tildeOmes} in Section 5.

\section{Proofs of Theorems \ref{gen-mes}, \ref{gen-sigma-mes}, \ref{gen-Nmes} and \ref{gen-Omes}}
 This section is devoted to giving the proofs of Theorems \ref{gen-mes}, \ref{gen-sigma-mes}, \ref{gen-Nmes} and \ref{gen-Omes}. Throughout this section, we assume that $r\geq 1$. We first present two proofs of Theorem \ref{gen-mes}.

{\noindent \bf The first proof of Theorem \ref{gen-mes}.}  For $k\geq 0$, $\pi$ is an overpartition with $mes_{r,A,a}(\pi)=kA+a$ if and only if
\[f_{\overline{iA+a}}(\pi)+f_{iA+a}(\pi)\geq r\quad\text{for}\quad 0\leq i\leq k-1,\]
and
\[f_{\overline{kA+a}}(\pi)+f_{kA+a}(\pi)\leq r-1.\]

More precisely, we obtain that $\pi$ is an overpartition with $mes_{r,A,a}(\pi)=kA+a$ if and only if
\[f_{\overline{iA+a}}(\pi)=0\text{ and }f_{iA+a}(\pi)\geq r,\text{ or } f_{\overline{iA+a}}(\pi)=1\text{ and }f_{iA+a}(\pi)\geq r-1,\text{ for } 0\leq i\leq k-1,\]
and
\[f_{\overline{kA+a}}(\pi)=0\text{ and }f_{kA+a}(\pi)\leq r-1,\text{ or } f_{\overline{kA+a}}(\pi)=1\text{ and }f_{kA+a}(\pi)\leq r-2.\]

So, we get
\begin{align*}
&\sum_{\pi\in\mathcal{\overline{P}}}z^{mes_{r,A,a}(\pi)}q^{|\pi|}\\
&=\sum_{k=0}^{\infty}z^{kA+a}2^kq^{r[a+(A+a)+\cdots+(k-1)A+a]}(1+2q^{kA+a}+\cdots+2q^{(r-1)(kA+a)})\\
&\qquad\times\frac{1-q^{kA+a}}{(q;q)_{\infty}}\frac{(-q;q)_{\infty}}{(-q^a;q^A)_{k+1}}\\
&=\frac{(-q;q)_{\infty}}{(q;q)_{\infty}}\sum_{k=0}^{\infty}z^{kA+a}\frac{2^kq^{r[A\binom{k}{2}+ka]}}{(-q^a;q^A)_{k+1}}\left(1+\frac{2q^{kA+a}-2q^{r(kA+a)}}{1-q^{kA+a}}\right)(1-q^{kA+a})\\
&=\frac{(-q;q)_{\infty}}{(q;q)_{\infty}}\sum_{k=0}^{\infty}z^{kA+a}\frac{2^kq^{r[A\binom{k}{2}+ka]}}{(-q^a;q^A)_{k+1}}(1+q^{kA+a}-2q^{r(kA+a)})\\
&=\frac{(-q;q)_{\infty}}{(q;q)_{\infty}}\sum_{k=0}^{\infty}z^{kA+a}\bigg[\frac{2^kq^{r[A\binom{k}{2}+ka]}}{(-q^a;q^A)_k}-\frac{2^{k+1}q^{r[A\binom{k+1}{2}+(k+1)a]}}{(-q^a;q^A)_{k+1}}\bigg].
\end{align*}

The proof is complete.  \qed

{\noindent \bf The second proof of Theorem \ref{gen-mes}.} For $k\geq 0$, $\pi$ is an overpartition with $mes_{r,A,a}(\pi)\geq kA+a$ if and only if
\[f_{\overline{iA+a}}(\pi)+f_{iA+a}(\pi)\geq r\quad\text{for}\quad 0\leq i\leq k-1.\]
This implies that the generating function for the overpartitions $\pi$ with $mes_{r,A,a}(\pi)\geq kA+a$ is
\[\frac{(-q;q)_{\infty}}{(q;q)_{\infty}}\frac{2^kq^{r[A\binom{k}{2}+ka]}}{(-q^a;q^A)_k},\]
and so the generating function for the overpartitions $\pi$ with $mes_{r,A,a}(\pi)\geq (k+1)A+a$ is
\[\frac{(-q;q)_{\infty}}{(q;q)_{\infty}}\frac{2^{k+1}q^{r[A\binom{k+1}{2}+(k+1)a]}}{(-q^a;q^A)_{k+1}}.\]
Hence, the generating function for the overpartitions $\pi$ with $mes_{r,A,a}(\pi)=kA+a$ is
\[\frac{(-q;q)_{\infty}}{(q;q)_{\infty}}\frac{2^kq^{r[A\binom{k}{2}+ka]}}{(-q^a;q^A)_k}-\frac{(-q;q)_{\infty}}{(q;q)_{\infty}}\frac{2^{k+1}q^{r[A\binom{k+1}{2}+(k+1)a]}}{(-q^a;q^A)_{k+1}}.\]
This completes the proof.  \qed

Then, we give two proofs of Theorem \ref{gen-sigma-mes}.

{\noindent \bf The first proof of Theorem \ref{gen-sigma-mes}.} Clearly, we have
\[\sum_{n=0}^\infty\sigma mes_{r,A,a}(n)q^n=\frac{\partial}{\partial z}\bigg|_{z=1}\left(\sum_{\pi\in\mathcal{\overline{P}}}z^{mes_{r,A,a}(\pi)}q^{|\pi|}\right).\]

By Theorem \ref{gen-mes}, we get
\begin{align*}
&\sum_{n=0}^\infty\sigma mes_{r,A,a}(n)q^n\\
&=\frac{(-q;q)_{\infty}}{(q;q)_{\infty}}\sum_{k=0}^{\infty}(kA+a)\bigg[\frac{2^kq^{r[A\binom{k}{2}+ka]}}{(-q^a;q^A)_k}-\frac{2^{k+1}
q^{r[A\binom{k+1}{2}+(k+1)a]}}{(-q^a;q^A)_{k+1}}\bigg]\nonumber\\
&=\frac{(-q;q)_{\infty}}{(q;q)_{\infty}}\bigg[\sum_{k=0}^{\infty}(kA+a)\frac{2^kq^{r[A\binom{k}{2}+ka]}}{(-q^a;q^A)_k}-\sum_{k=0}^{\infty}(kA+a)\frac{2^{k+1}
q^{r[A\binom{k+1}{2}+(k+1)a]}}{(-q^a;q^A)_{k+1}}\bigg]\nonumber\\
&=\frac{(-q;q)_{\infty}}{(q;q)_{\infty}}\bigg[\sum_{k=0}^{\infty}(kA+a)\frac{2^kq^{r[A\binom{k}{2}+ka]}}{(-q^a;q^A)_k}-\sum_{k=1}^{\infty}((k-1)A+a)
\frac{2^kq^{r[A\binom{k}{2}+ka]}}{(-q^a;q^A)_k}\bigg]\nonumber\\
&=\frac{(-q;q)_{\infty}}{(q;q)_{\infty}}\bigg[a+A\sum_{k=1}^{\infty}\frac{2^kq^{r[A\binom{k}{2}+ka]}}{(-q^a;q^A)_k}\bigg].
\end{align*}

The proof is complete.  \qed

{\noindent \bf The second proof of Theorem \ref{gen-sigma-mes}.} For $k\geq 0$ and $n\geq 0$, let $M_{r,A,a}(k,n)$ be the number of overpartitions $\pi$ of $n$ with  $mes_{r,A,a}(\pi)\geq kA+a$. In the second proof of Theorem \ref{gen-mes}, we have shown that for $k\geq 0$,
\begin{equation}\label{M-gen}
\sum_{n=0}^{\infty}M_{r,A,a}(k,n)q^n=\frac{(-q;q)_{\infty}}{(q;q)_{\infty}}\frac{2^{k}q^{r[A\binom{k}{2}+ka]}}{(-q^a;q^A)_{k}}.
\end{equation}
We proceed to show that for $n\geq0$,
\begin{equation}\label{relation-sigma-M}
\sigma mes_{r,A,a}(n)=A\sum_{k=1}^{\infty}M_{r,A,a}(k,n)+a\overline{p}(n).
\end{equation}

It is from the definition of $M_{r,A,a}(k,n)$ that for $k\geq 0$ and $n\geq 0$,
\[\sum_{\pi\in\overline{\mathcal{P}}(n),\ mes_{r,A,a}(\pi)=kA+a}1=M_{r,A,a}(k,n)-M_{r,A,a}(k+1,n).\]

Then, we have
\begin{align*}
\sigma mes_{r,A,a}(n)&=\sum_{\pi\in\overline{\mathcal{P}}(n)}mes_{r,A,a}(\pi)\nonumber\\
&=\sum_{k=0}^{\infty}\sum_{\pi\in\overline{\mathcal{P}}(n),\ mes_{r,A,a}(\pi)=kA+a}kA+a\\
&=A\sum_{k=0}^{\infty}k\sum_{\pi\in\overline{\mathcal{P}}(n),\ mes_{r,A,a}(\pi)=kA+a}1+a\sum_{k=0}^{\infty}\sum_{\pi\in\overline{\mathcal{P}}(n),\ mes_{r,A,a}(\pi)=kA+a}1\\
&=A\sum_{k=0}^{\infty}k \bigg[M_{r,A,a}(k,n)-M_{r,A,a}(k+1,n)\bigg]+a\overline{p}(n)\\
&=A\bigg[\sum_{k=1}^{\infty}k M_{r,A,a}(k,n)-\sum_{k=1}^{\infty}(k-1) M_{r,A,a}(k,n)\bigg]+a\overline{p}(n)\\
&=A\sum_{k=1}^{\infty}M_{r,A,a}(k,n)+a\overline{p}(n).
\end{align*}
We arrive at \eqref{relation-sigma-M}. Combining \eqref{gen-over-p}, \eqref{M-gen} and \eqref{relation-sigma-M}, we get
\begin{align*}
\sum_{n=0}^\infty\sigma mes_{r,A,a}(n)q^n&=A\sum_{n=0}^\infty \sum_{k=1}^{\infty}M_{r,A,a}(k,n)q^n+a\sum_{n=0}^\infty\overline{p}(n)q^n\\
&=A\sum_{k=1}^\infty \sum_{n=0}^{\infty}M_{r,A,a}(k,n)q^n+a\frac{(-q;q)_{\infty}}{(q;q)_{\infty}}\\
&=A\sum_{k=1}^\infty\frac{(-q;q)_{\infty}}{(q;q)_{\infty}}\frac{2^{k}q^{r[A\binom{k}{2}+ka]}}{(-q^a;q^A)_{k}}+a\frac{(-q;q)_{\infty}}{(q;q)_{\infty}}\\
&=\frac{(-q;q)_{\infty}}{(q;q)_{\infty}}\bigg[a+A\sum_{k=1}^{\infty}\frac{2^kq^{r[A\binom{k}{2}+ka]}}{(-q^a;q^A)_k}\bigg].
\end{align*}
This completes the proof.  \qed

Then, we will give two proofs of Theorem \ref{gen-Nmes} by using the similar arguments in the proofs of Theorem \ref{gen-mes}.

{\noindent \bf The first proof of Theorem \ref{gen-Nmes}.}  In view of the definition of ${Nmes}_{r,A,a}$, we see that for $k\geq 0$, $\pi$ is an overpartition  in ${Nmes}_{r,A,a}$ with $mes_{r,A,a}(\pi)=kA+a$ if and only if
\[f_{\overline{iA+a}}(\pi)=0\text{ and }f_{iA+a}(\pi)\geq r, \text{ for } 0\leq i\leq k-1,\]
and
\[f_{\overline{kA+a}}(\pi)=0\text{ and }f_{kA+a}(\pi)\leq r-1,\text{ or } f_{\overline{kA+a}}(\pi)=1\text{ and }f_{kA+a}(\pi)\leq r-2.\]
So, we get
\begin{align*}
&\sum_{\pi\in{Nmes}_{r,A,a}}z^{mes_{r,A,a}(\pi)}q^{|\pi|}\\
&=\sum_{k=0}^{\infty}z^{kA+a}q^{r[a+(A+a)+\cdots+(k-1)A+a]}(1+2q^{kA+a}+\cdots+2q^{(r-1)(kA+a)})\\
&\qquad\times\frac{1-q^{kA+a}}{(q;q)_{\infty}}\frac{(-q;q)_{\infty}}{(-q^a;q^A)_{k+1}}\\
&=\frac{(-q;q)_{\infty}}{(q;q)_{\infty}}\sum_{k=0}^{\infty}z^{kA+a}\frac{q^{r[A\binom{k}{2}+ka]}}{(-q^a;q^A)_{k+1}}(1+q^{kA+a}-2q^{r(kA+a)})\\
&=\frac{(-q;q)_{\infty}}{(q;q)_{\infty}}\sum_{k=0}^{\infty}z^{kA+a}\bigg[\frac{q^{r[A\binom{k}{2}+ka]}}{(-q^a;q^A)_k}-\frac{2q^{r[A\binom{k+1}{2}+(k+1)a]}}{(-q^a;q^A)_{k+1}}\bigg].
\end{align*}

The proof is complete.  \qed

{\noindent \bf The second proof of Theorem \ref{gen-Nmes}.} For $k\geq 0$, it is clear that $\pi$ is an overpartition such that $mes_{r,A,a}(\pi)\geq kA+a$ and all parts $\equiv a\pmod A$ of size less than $kA+a$ in $\pi$ are non-overlined if and only if
\[f_{\overline{iA+a}}(\pi)=0\text{ and }f_{iA+a}(\pi)\geq r, \text{ for } 0\leq i\leq k-1,\]
and so the generating function for such overpartitions $\pi$ is
\[\frac{(-q;q)_{\infty}}{(q;q)_{\infty}}\frac{q^{r[A\binom{k}{2}+ka]}}{(-q^a;q^A)_k}.\]

Note that $\pi$ is an overpartition such that $mes_{r,A,a}(\pi)\geq (k+1)A+a$ and all parts $\equiv a\pmod A$ of size less than $kA+a$ in $\pi$ are non-overlined if and only if
\[f_{\overline{iA+a}}(\pi)=0\text{ and }f_{iA+a}(\pi)\geq r, \text{ for } 0\leq i\leq k-1,\]
and
\[f_{\overline{kA+a}}(\pi)=0\text{ and }f_{kA+a}(\pi)\geq r,\text{ or } f_{\overline{kA+a}}(\pi)=1\text{ and }f_{kA+a}(\pi)\geq r-1,\]
so we see that the generating function for such overpartitions $\pi$ is
\[\frac{(-q;q)_{\infty}}{(q;q)_{\infty}}\frac{2q^{r[A\binom{k+1}{2}+(k+1)a]}}{(-q^a;q^A)_{k+1}}.\]
Hence, the generating function for the overpartitions $\pi$ such that $mes_{r,A,a}(\pi)=kA+a$ and all parts $\equiv a\pmod A$ of size less than $kA+a$ in $\pi$ are non-overlined is
\[\frac{(-q;q)_{\infty}}{(q;q)_{\infty}}\frac{q^{r[A\binom{k}{2}+ka]}}{(-q^a;q^A)_k}-\frac{(-q;q)_{\infty}}{(q;q)_{\infty}}\frac{2q^{r[A\binom{k+1}{2}+(k+1)a]}}{(-q^a;q^A)_{k+1}}.\]
This completes the proof.  \qed

We conclude this section with two proofs of Theorem \ref{gen-Omes}.

{\noindent \bf The first proof of Theorem \ref{gen-Omes}.}  By definition, we find that for $k\geq 0$, $\pi$ is an overpartition in ${Omes}_{1,A,a}$ with $mes_{1,A,a}(\pi)=kA+a$ if and only if
\[f_{\overline{iA+a}}(\pi)=1\text{ and }f_{iA+a}(\pi)=0, \text{ for } 0\leq i\leq k-1,\]
and
\[f_{\overline{kA+a}}(\pi)=0\text{ and }f_{kA+a}(\pi)=0.\]
So, we get
\begin{align*}
&\sum_{\pi\in Omes_{1,A,a}}z^{mes_{1,A,a}(\pi)}q^{|\pi|}\\
&=\sum_{k=0}^{\infty}z^{kA+a}q^{[a+(A+a)+\cdots+(k-1)A+a]}
\frac{(q^a;q^A)_{k+1}}{(q;q)_{\infty}}\frac{(-q;q)_{\infty}}{(-q^a;q^A)_{k+1}}\nonumber\\
&=\frac{(-q;q)_{\infty}}{(q;q)_{\infty}}\sum_{k=0}^{\infty}z^{kA+a}q^{A\binom{k}{2}+ka} \frac{(q^a;q^A)_{k}(1+q^{kA+a}-2q^{kA+a})}{(-q^a;q^A)_{k+1}}\\
&=\frac{(-q;q)_{\infty}}{(q;q)_{\infty}}\sum_{k=0}^{\infty}z^{kA+a}\bigg[q^{A\binom{k}{2}+ka} \frac{(q^a;q^A)_k}{(-q^a;q^A)_k}-
2q^{A\binom{k+1}{2}+(k+1)a} \frac{(q^a;q^A)_k}{(-q^a;q^A)_{k+1}}\bigg].
\end{align*}

The proof is complete.  \qed

{\noindent \bf The second proof of Theorem \ref{gen-Omes}.} For $k\geq 0$, it is clear that $\pi$ is an overpartition with $mes_{1,A,a}(\pi)\geq kA+a$ and all parts $\equiv a\pmod A$ of size less than $kA+a$ in $\pi$ are overlined if and only if
\[f_{\overline{iA+a}}(\pi)=1\text{ and }f_{iA+a}(\pi)=0, \text{ for } 0\leq i\leq k-1,\]
and so the generating function for such overpartitions $\pi$ is
\[\frac{(-q;q)_{\infty}}{(q;q)_{\infty}}q^{A\binom{k}{2}+ka} \frac{(q^a;q^A)_k}{(-q^a;q^A)_k}.\]

Note that $\pi$ is an overpartition with $mes_{1,A,a}(\pi)\geq (k+1)A+a$ and all parts $\equiv a\pmod A$ of size less than $kA+a$ in $\pi$ are overlined if and only if
\[f_{\overline{iA+a}}(\pi)=1\text{ and }f_{iA+a}(\pi)=0, \text{ for } 0\leq i\leq k-1,\]
and
\[f_{\overline{kA+a}}(\pi)=0\text{ and }f_{kA+a}(\pi)\geq 1,\text{ or } f_{\overline{kA+a}}(\pi)=1\text{ and }f_{kA+a}(\pi)\geq 0,\]
we see that the generating function for such overpartitions $\pi$ is
\[\frac{(-q;q)_{\infty}}{(q;q)_{\infty}}2q^{A\binom{k+1}{2}+(k+1)a} \frac{(q^a;q^A)_k}{(-q^a;q^A)_{k+1}}.\]
Hence, the generating function for the overpartitions $\pi$ such that ${mes}_{1,A,a}(\pi)=kA+a$ and all parts $\equiv a\pmod A$ of size less than $kA+a$ in $\pi$ are overlined is
\[\frac{(-q;q)_{\infty}}{(q;q)_{\infty}}q^{A\binom{k}{2}+ka} \frac{(q^a;q^A)_k}{(-q^a;q^A)_k}-\frac{(-q;q)_{\infty}}{(q;q)_{\infty}}2q^{A\binom{k+1}{2}+(k+1)a} \frac{(q^a;q^A)_k}{(-q^a;q^A)_{k+1}}.\]
This completes the proof.  \qed

\section{Proofs of Theorems \ref{gen-overmes}, \ref{gen-sigma-overmes} and \ref{gen-overNmes}}

In this section, we aim to show Theorems \ref{gen-overmes}, \ref{gen-sigma-overmes} and \ref{gen-overNmes}. We first prove Theorem \ref{gen-overmes} with the similar arguments in the proofs of Theorem \ref{gen-mes}. In the remaining of this article, we fix $r\geq 2$.

{\noindent \bf The first proof of Theorem \ref{gen-overmes}.}  For $k\geq0$, $\pi$ is an overpartition with $\overline{mes}_{r,A,a}(\pi)=kA+a$ if and only if
\[f_{\overline{iA+a}}(\pi)=1\text{ and }f_{iA+a}(\pi)\geq r-1,\text{ for } 0\leq i\leq k-1,\]
and
\[f_{\overline{kA+a}}(\pi)=0\text{ or }f_{kA+a}(\pi)\leq r-2.\]
Note that $f_{\overline{kA+a}}(\pi)\leq 1$, we see that if $f_{\overline{kA+a}}(\pi)+f_{kA+a}(\pi)\geq r$, then $f_{kA+a}(\pi)\geq r-1$, and so $f_{\overline{kA+a}}(\pi)=0$.  Hence, we have

\begin{align*}
&\sum_{\pi\in\mathcal{\overline{P}}}z^{\overline{mes}_{r,A,a}(\pi)}q^{|\pi|}\\
&=\sum_{k=0}^{\infty}z^{kA+a}q^{r[a+(A+a)+\cdots+(k-1)A+a]}\left(1+2q^{kA+a}+\cdots+2q^{(r-1)(kA+a)}+
\sum_{i=r}^{\infty}q^{i(kA+a)}\right)\\
&\qquad\qquad\times\frac{1-q^{kA+a}}{(q;q)_{\infty}}\frac{(-q;q)_{\infty}}{(-q^a;q^A)_{k+1}}\nonumber\\
&=\frac{(-q;q)_{\infty}}{(q;q)_{\infty}}\sum_{k=0}^{\infty}z^{kA+a}\frac{q^{r[A\binom{k}{2}+ka]}}{(-q^a;q^A)_{k+1}}(1+q^{kA+a}-q^{r(kA+a)})\\
&=\frac{(-q;q)_{\infty}}{(q;q)_{\infty}}\sum_{k=0}^{\infty}z^{kA+a}\bigg[\frac{q^{r[A\binom{k}{2}+ka]}}{(-q^a;q^A)_k}-\frac{q^{r[A\binom{k+1}{2}+(k+1)a]}}{(-q^a;q^A)_{k+1}}\bigg].
\end{align*}

The proof is complete.  \qed

{\noindent \bf The second proof of Theorem \ref{gen-overmes}.} For $k\geq 0$, $\pi$ is an overpartition  with $\overline{mes}_{r,A,a}(\pi)\geq kA+a$ if and only if
\[f_{\overline{iA+a}}(\pi)=1\text{ and }f_{iA+a}(\pi)\geq r-1,\text{ for } 0\leq i\leq k-1.\]
Then, the generating function for the overpartitions $\pi$ with $\overline{mes}_{r,A,a}(\pi)\geq kA+a$ is
\[\frac{(-q;q)_{\infty}}{(q;q)_{\infty}}\frac{q^{r[A\binom{k}{2}+ka]}}{(-q^a;q^A)_k},\]
and so the generating function for the overpartitions $\pi$ with $\overline{mes}_{r,A,a}(\pi)\geq (k+1)A+a$ is
\[\frac{(-q;q)_{\infty}}{(q;q)_{\infty}}\frac{q^{r[A\binom{k+1}{2}+(k+1)a]}}{(-q^a;q^A)_{k+1}}.\]
Hence, the generating function for the overpartitions $\pi$ with $\overline{mes}_{r,A,a}(\pi)=kA+a$ is
\[\frac{(-q;q)_{\infty}}{(q;q)_{\infty}}\frac{q^{r[A\binom{k}{2}+ka]}}{(-q^a;q^A)_k}-\frac{(-q;q)_{\infty}}{(q;q)_{\infty}}\frac{q^{r[A\binom{k+1}{2}+(k+1)a]}}{(-q^a;q^A)_{k+1}}.\]
This completes the proof.  \qed

Then, we give two proofs of Theorem \ref{gen-sigma-overmes} by using the similar arguments in the proofs of Theorem \ref{gen-sigma-mes}.

{\noindent \bf The first proof of Theorem \ref{gen-sigma-overmes}.} Appealing Theorem \ref{gen-overmes}, we have
\begin{align*}
\sum_{n=0}^\infty\sigma \overline{mes}_{r,A,a}(n)q^n&=\frac{\partial}{\partial z}\bigg|_{z=1}\left(\sum_{\pi\in\mathcal{\overline{P}}}z^{\overline{mes}_{r,A,a}(\pi)}q^{|\pi|}\right)\\
&=\frac{{(-q;q)_{\infty}}}{{(q;q)_{\infty}}}\sum_{k=0}^{\infty}(kA+a)\bigg[\frac{q^{r[A\binom{k}{2}+ka]}}{(-q^{a};q^{A})_{k}}-\frac{q^{r[A\binom{k+1}{2}+(k+1)a]}}{(-q^{a};q^{A})_{k+1}}\bigg]\\
&=\frac{{(-q;q)_{\infty}}}{{(q;q)_{\infty}}}\bigg[\sum_{k=0}^{\infty}(kA+a)\frac{q^{r[A\binom{k}{2}+ka]}}{(-q^{a};q^{A})_{k}}-\sum_{k=1}^{\infty}((k-1)A+a)\frac{q^{r[A\binom{k}{2}+ka]}}{(-q^{a};q^{A})_{k}}\bigg]\\
&=\frac{{(-q;q)_{\infty}}}{{(q;q)_{\infty}}}\bigg[a+A\sum_{k=1}^{\infty}\frac{q^{r[A\binom{k}{2}+ka]}}{(-q^{a};q^{A})_{k}}\bigg].
\end{align*}

The proof is complete.  \qed

{\noindent \bf The second proof of Theorem \ref{gen-sigma-overmes}.} For $k\geq 0$ and $n\geq 0$, let $\overline{M}_{r,A,a}(k,n)$ be the number of overpartitions $\pi$ of $n$ with  $\overline{mes}_{r,A,a}(\pi)\geq kA+a$. In the second proof of Theorem \ref{gen-overmes}, we have shown that for $k\geq 0$,
\begin{equation}\label{overM-gen}
\sum_{n=0}^\infty\overline{M}_{r,A,a}(k,n)q^{n}=\frac{{(-q;q)_{\infty}}}{{(q;q)_{\infty}}}\frac{q^{r[A\binom{k}{2}+ka]}}{(-q^{a};q^{A})_{k}}.
\end{equation}

With a similar argument in the proof of \eqref{relation-sigma-M}, we can obtain that for $n\geq 0$,
\begin{equation}\label{relation-sigma-overM}
\sigma\overline{mes}_{r,A,a}(n)=A\sum_{k=1}^{\infty}\overline{M}_{r,A,a}(k,n)+a\overline{p}(n).
\end{equation}

 Combining \eqref{gen-over-p}, \eqref{overM-gen} and \eqref{relation-sigma-overM}, we get
\begin{align*}
\sum_{n=0}^\infty\sigma \overline{mes}_{r,A,a}(n)q^n&=A\sum_{n=0}^\infty \sum_{k=1}^{\infty}\overline{M}_{r,A,a}(k,n)q^n+a\sum_{n=0}^\infty\overline{p}(n)q^n\\
&=A\sum_{k=1}^\infty \sum_{n=0}^{\infty}\overline{M}_{r,A,a}(k,n)q^n+a\frac{(-q;q)_{\infty}}{(q;q)_{\infty}}\\
&=A\sum_{k=1}^\infty\frac{{(-q;q)_{\infty}}}{{(q;q)_{\infty}}}\frac{q^{r[A\binom{k}{2}+ka]}}{(-q^{a};q^{A})_{k}}+a\frac{(-q;q)_{\infty}}{(q;q)_{\infty}}\\
&=\frac{{(-q;q)_{\infty}}}{{(q;q)_{\infty}}}\bigg[a+A\sum_{k=1}^{\infty}\frac{q^{r[A\binom{k}{2}+ka]}}{(-q^{a};q^{A})_{k}}\bigg].
\end{align*}
This completes the proof.  \qed

We conclude this section with two proofs of Theorem \ref{gen-overNmes}.

{\noindent \bf The first proof of Theorem \ref{gen-overNmes}.}  Under the definition of ${N\overline{mes}}_{r,A,a}$, we find that for $k\geq 0$, $\pi$ is an overpartition in ${N\overline{mes}}_{r,A,a}$ with  $\overline{mes}_{r,A,a}(\pi)=kA+a$ if and only if
\[f_{\overline{iA+a}}(\pi)=1\text{ and }f_{iA+a}(\pi)\geq r-1,\text{ for } 0\leq i\leq k-1,\]
and
\[f_{kA+a}(\pi)\leq r-2.\]
So, we get
\begin{align*}
\sum_{\pi\in N\overline{mes}_{r,A,a}}z^{\overline{mes}_{r,A,a}(\pi)}q^{|\pi|}&=\sum_{k=0}^{\infty}z^{kA+a}q^{r[a+(A+a)+\cdots+(k-1)A+a]}(1+q^{kA+a}+\cdots+q^{(r-2)(kA+a)})\\
&\qquad\qquad\times\frac{1-q^{kA+a}}{{(q;q)_{\infty}}}\frac{{(-q;q)_{\infty}}}{(-q^{a};q^{A})_{k}}\\
&=\frac{{(-q;q)_{\infty}}}{{(q;q)_{\infty}}}\sum_{k=0}^{\infty}z^{kA+a}\frac{q^{r[A\binom{k}{2}+ka]}}{(-q^{a};q^{A})_{k}}(1-q^{(r-1)(kA+a)})\\
&=\frac{{(-q;q)_{\infty}}}{{(q;q)_{\infty}}}\sum_{k=0}^{\infty}z^{kA+a}\bigg[\frac{q^{r[A\binom{k}{2}+ka]}}{(-q^{a};q^{A})_{k}}-\frac{q^{r[A\binom{k}{2}+ka]+(r-1)(kA+a)}}{(-q^{a};q^{A})_{k}}\bigg].
\end{align*}

The proof is complete.  \qed

{\noindent \bf The second proof of Theorem \ref{gen-overNmes}.} For $k\geq 0$, in the second proof of Theorem \ref{gen-overmes}, we have proved that
\[\frac{(-q;q)_{\infty}}{(q;q)_{\infty}}\frac{q^{r[A\binom{k}{2}+ka]}}{(-q^a;q^A)_k}\]
is the generating function for the overpartitions $\pi$ with $\overline{mes}_{r,A,a}(\pi)\geq kA+a$. Moveover, we see that
\[\frac{{(-q;q)_{\infty}}}{{(q;q)_{\infty}}}\frac{q^{r[A\binom{k}{2}+ka]+(r-1)(kA+a)}}{(-q^{a};q^{A})_{k}}\]
is the generating function for the overpartitions  $\pi$  such that $\overline{mes}_{r,A,a}(\pi)\geq kA+a$ and
  there are  at  least  $r-1$  non-overlined  parts  of  size  $kA+a$ in $\pi$. So, we see that
  \begin{equation}\label{overN-substract}
  \frac{(-q;q)_{\infty}}{(q;q)_{\infty}}\frac{q^{r[A\binom{k}{2}+ka]}}{(-q^a;q^A)_k}-\frac{{(-q;q)_{\infty}}}{{(q;q)_{\infty}}}\frac{q^{r[A\binom{k}{2}+ka]+(r-1)(kA+a)}}{(-q^{a};q^{A})_{k}}
  \end{equation}
  is the generating function for the overpartitions  $\pi$  such that
  \begin{itemize}
  \item[(1)] $\overline{mes}_{r,A,a}(\pi)\geq kA+a$;
  \item[(2)] there are  no  $r-1$  non-overlined  parts  of  size  $kA+a$ in $\pi$.
  \end{itemize}
  By definition, we know that the conditions (1) and (2) yields $\overline{mes}_{r,A,a}(\pi)=kA+a$, and so \eqref{overN-substract} is the  the generating function for the overpartitions  $\pi$ in $N\overline{mes}_{r,A,a}$ with $\overline{mes}_{r,A,a}(\pi)=kA+a$. This completes the proof.  \qed

\section{Proofs of Theorems \ref{gen-tildemes}, \ref{gen-sigma-tildemes}, \ref{gen-tildeNmes} and \ref{gen-tildeOmes}}

The objective of this section is to present the proofs of Theorems \ref{gen-tildemes}, \ref{gen-sigma-tildemes}, \ref{gen-tildeNmes} and \ref{gen-tildeOmes}. Mimicking the proofs of Theorems \ref{gen-mes} \ref{gen-overmes}, we give two proofs of Theorem \ref{gen-tildemes}.

{\noindent \bf The first proof of Theorem \ref{gen-tildemes}.}  For $k\geq 0$, we find that $\pi$ is an overpartition with $\widetilde{mes}_{r,A,a}(\pi)=kA+a$ if and only if
\[f_{\overline{iA+a}}(\pi)=1\text{ or }f_{iA+a}(\pi)\geq r-1,\text{ for } 0\leq i\leq k-1,\]
and
\[f_{\overline{kA+a}}(\pi)=0\text{ and }f_{kA+a}(\pi)\leq r-2.\]
Moreover, we see that for $0\leq i\leq k-1$, if $f_{\overline{iA+a}}(\pi)+f_{iA+a}(\pi)\leq r-2$ then we have $f_{\overline{iA+a}}(\pi)=1$. So, we get
\begin{align*}
	&\sum_{\pi \in \overline{P}}z^{\widetilde{mes}_{r,A,a}(\pi)} q^{|\pi|}\\
&=\sum_{k=0}^{\infty }
	z^{kA+a}\prod_{i=0}^{k-1}\bigg[q^{iA+a}+q^{2(iA+a)}+\cdots+q^{(r-2)(iA+a)} +\sum_{m=r-1}^{\infty } 2q^{m(iA+a)}\bigg]\\
	&\qquad\qquad\times(1+q^{kA+a}+\cdots+q^{(r-2)(kA+a)})\frac{(q^{a};q^{A})_{k+1}}{(q;q)_{\infty }} \frac{(-q;q)_{\infty }}{(-q^{a};q^{a})_{k+1}}\nonumber\\
	&=\frac{(-q;q)_{\infty }}{(q;q)_{\infty}}\sum_{k=0}^{\infty }z^{kA+a}\prod_{i=0}^{k-1} \frac{q^{iA+a}+q^{(r-1)(iA+a)}}{1-q^{iA+a}}
	\frac{1-q^{(r-1)(kA+a)}}{1-q^{kA+a}} \frac{(q^{a};q^{A})_{k+1}}{(-q^{a};q^{A})_{k+1}} \\
	&=\frac{(-q;q)_{\infty}}{(q;q)_{\infty}}\sum_{k=0}^{\infty }z^{kA+a}\prod_{i=0}^{k-1}\bigg[ q^{iA+a}(1+q^{(r-2)(iA+a)})\bigg]
	\frac{(1+q^{kA+a}-q^{kA+a}-q^{(r-1)(kA+a)})}{(-q^{a};q^{A})_{k+1}}\\
	&=\frac{(-q;q)_{\infty}}{(q;q)_{\infty}}\sum_{k=0}^{\infty}z^{kA+a}\bigg[\frac{q^{A[{k\choose 2}+ka ]}}{(-q^{a};q^{A})_{k}}(-q^{(r-2)a};q^{(r-2)A})_{k}\\
&\qquad\qquad\qquad\qquad\qquad-\frac{q^{A[{k+1\choose 2}+(k+1)a ]}}{(-q^{a};q^{A})_{k+1}}(-q^{(r-2)a};q^{(r-2)A})_{k+1}\bigg].
\end{align*}

The proof is complete.  \qed

{\noindent \bf The second proof of Theorem \ref{gen-tildemes}.} For $k\geq 0$,  we find that $\pi$ is an overpartition with $\widetilde{mes}_{r,A,a}(\pi)\geq kA+a$ if and only if
\[f_{\overline{iA+a}}(\pi)=1\text{ or }f_{iA+a}(\pi)\geq r-1,\text{ for } 0\leq i\leq k-1.\]

It is clear that
\[q^{A[{k\choose 2}+ka]}(-q^{(r-2)a};q^{(r-2)A})_{k}\]
is the generating function for the overpartitions $\pi$ such that all parts of $\pi$ are congruent to $a$ modulo $A$, all parts of $\pi$ are of size not exceeding $(k-1)A+a$, and
\[f_{\overline{iA+a}}(\pi)=1\text{ and }f_{iA+a}(\pi)=0,\text{ or }f_{\overline{iA+a}}(\pi)=0\text{ and }f_{iA+a}(\pi)=r-1,\text{ for } 0\leq i\leq k-1.\]

This implies that the generating function for the overpartitions $\pi$ with $\widetilde{mes}_{r,A,a}(\pi)\geq kA+a$ is
\[\frac{(-q;q)_{\infty }}{(q;q)_{\infty }}\frac{q^{A[{k\choose 2}+ka]}}{(-q^{a};q^{A})_{k}} (-q^{(r-2)a};q^{(r-2)A})_{k},\]
and so the generating function for the overpartitions $\pi$ with $\widetilde{mes}_{r,A,a}(\pi)\geq (k+1)A+a$ is
\[\frac{(-q;q)_{\infty }}{(q;q)_{\infty }}\frac{q^{A[{k+1\choose 2}+(k+1)a]}}{(-q^{a};q^{A})_{k+1}} (-q^{(r-2)a};q^{(r-2)A})_{k+1}.\]
Hence, the generating function for the overpartitions $\pi$ with $\widetilde{mes}_{r,A,a}(\pi)=kA+a$ is
\[\frac{(-q;q)_{\infty }}{(q;q)_{\infty }}\frac{q^{A[{k\choose 2}+ka]}}{(-q^{a};q^{A})_{k}} (-q^{(r-2)a};q^{(r-2)A})_{k}-\frac{(-q;q)_{\infty }}{(q;q)_{\infty }}\frac{q^{A[{k+1\choose 2}+(k+1)a]}}{(-q^{a};q^{A})_{k+1}} (-q^{(r-2)a};q^{(r-2)A})_{k+1}.\]
This completes the proof.  \qed

Then, we give two proofs of Theorem \ref{gen-sigma-tildemes} by using the similar arguments in the proofs of Theorems \ref{gen-sigma-mes} and \ref{gen-sigma-overmes}.

{\noindent \bf The first proof of Theorem \ref{gen-sigma-tildemes}.} In view of Theorem \ref{gen-tildemes}, we have
\begin{align*}
\sum_{n=0}^\infty\sigma \widetilde{mes}_{r,A,a}(n)q^n&=\frac{\partial}{\partial z}\bigg|_{z=1}\left(\sum_{\pi\in\mathcal{\overline{P}}}z^{\widetilde{mes}_{r,A,a}(\pi)}q^{|\pi|}\right)\\
&=\frac{(-q;q)_{\infty }}{(q;q)_{\infty }} \bigg[\sum_{k=0}^{\infty } (kA+a)\frac{q^{A[{k\choose 2}+ka]}}{(-q^{a};q^{A})_{k}}
	(-q^{(r-2)a};q^{(r-2)A})_{k}\\
	&\qquad\qquad\qquad-\sum_{k=1}^{\infty }((k-1)A+a)\frac{q^{A[{k\choose 2}+ka]}}{(-q^{a};q^{A})_{k}}(-q^{(r-2)a};q^{(r-2)A})_{k}\bigg] \\
	=&\frac{(-q;q)_{\infty }}{(q;q)_{\infty }} \bigg[a+A\sum_{k=1}^{\infty }\frac{q^{A[{k\choose 2}+ka]}}{(-q^{a};q^{A})_{k}}(-q^{(r-2)a};q^{(r-2)A})_{k}\bigg].
\end{align*}
The proof is complete.  \qed

{\noindent \bf The second proof of Theorem \ref{gen-sigma-tildemes}.} For $k\geq 0$ and $n\geq0$, let $\widetilde{M}_{r,A,a}(k,n)$ be the number of overpartitions $\pi$ of $n$ with  $\widetilde{mes}_{r,A,a}(\pi)\geq kA+a$. In the second proof of Theorem \ref{gen-tildemes}, we have shown that for $k\geq 0$,
\begin{equation}\label{tildeM-gen}
\sum_{n=0}^\infty\widetilde{M}_{r,A,a}(k,n)q^{n}=\frac{(-q;q)_{\infty}}{(q;q)_{\infty}}\frac{q^{A[{k\choose 2}+ka]}}{(-q^{a};q^{A})_{k}} (-q^{(r-2)a};q^{(r-2)A})_{k}.
\end{equation}

With a similar argument in the proof of \eqref{relation-sigma-M}, we can get that for $n\geq 0$,
\begin{equation}\label{relation-sigma-tildeM}
\sigma\widetilde{mes}_{r,A,a}(n)=A\sum_{k=1}^{\infty}\widetilde{M}_{r,A,a}(k,n)+a\overline{p}(n).
\end{equation}

It follows from \eqref{gen-over-p}, \eqref{tildeM-gen} and \eqref{relation-sigma-tildeM} that
\begin{align*}
\sum_{n=0}^\infty\sigma \widetilde{mes}_{r,A,a}(n)q^n&=A\sum_{n=0}^\infty \sum_{k=1}^{\infty}\widetilde{M}_{r,A,a}(k,n)q^n+a\sum_{n=0}^\infty\overline{p}(n)q^n\\
&=A\sum_{k=1}^\infty \sum_{n=0}^{\infty}\widetilde{M}_{r,A,a}(k,n)q^n+a\frac{(-q;q)_{\infty}}{(q;q)_{\infty}}\\
&=A\sum_{k=1}^\infty\frac{(-q;q)_{\infty}}{(q;q)_{\infty}}\frac{q^{A[{k\choose 2}+ka]}}{(-q^{a};q^{A})_{k}} (-q^{(r-2)a};q^{(r-2)A})_{k}+a\frac{(-q;q)_{\infty}}{(q;q)_{\infty}}\\
&=\frac{{(-q;q)_{\infty}}}{{(q;q)_{\infty}}}\bigg[a+A\sum_{k=1}^{\infty}\frac{q^{A[{k\choose 2}+ka]}}{(-q^{a};q^{A})_{k}}(-q^{(r-2)a};q^{(r-2)A})_{k}\bigg].
\end{align*}
This completes the proof.  \qed

Now, we are in a position to present two proofs of Theorem \ref{gen-tildeNmes}.

{\noindent \bf The first proof of Theorem \ref{gen-tildeNmes}.}  We find that for $k\geq 0$, $\pi$ is
 an overpartition in ${N\widetilde{mes}}_{r,A,a}$ with $\widetilde{mes}_{r,A,a}(\pi)=kA+a$ if and only if
\[f_{\overline{iA+a}}(\pi)=0\text{ and }f_{iA+a}(\pi)\geq r-1,\text{ for } 0\leq i\leq k-1,\]
and
\[f_{\overline{kA+a}}(\pi)=0\text{ and }f_{kA+a}(\pi)\leq r-2.\]
So, we get
	\begin{align*}
	&\sum_{\pi \in N\widetilde{mes}_{r,A,a} }z^{\widetilde{mes}_{r,A,a}(\pi)}q^{\left | \pi  \right | }\\
	&=\sum_{k=0}^{\infty} z^{kA+a}q^{(r-1)[a+(A+a)+\cdots+(k-1)A+a]}(1+q^{kA+a}+\cdots+q^{(r-2)(kA+a)})
	\frac{1-q^{kA+a}}{(q;q)_{\infty }} \frac{(-q;q)_{\infty }}{(-q^{a};q^{A})_{k+1 }} \\
    &=\frac{(-q;q)_{\infty }}{(q;q)_{\infty}} \sum_{k=0}^{\infty}z^{kA+a}\frac{q^{(r-1)[A{k\choose 2}+ka]}}{(-q^{a};q^{A})_{k+1 }}
	(1-q^{(r-1)(kA+a)})\\
	&=\frac{(-q;q)_{\infty }}{(q;q)_{\infty}}\sum_{k=0}^{\infty }z^{kA+a}\bigg[\frac{q^{(r-1)[A{k\choose 2}+ka]}}{(-q^{a};q^{A})_{k+1}}-\frac{q^{(r-1)[A{k+1\choose 2}+(k+1)a]}}{(-q^{a};q^{A})_{k+1 }}\bigg]. \nonumber
	\end{align*}
The proof is complete.  \qed

{\noindent \bf The second proof of Theorem \ref{gen-tildeNmes}.} For $k\geq 0$, it is clear that
\[\frac{1}{(q;q)_{\infty}}\frac{(-q;q)_{\infty}}{(-q^{a};q^{A})_{k+1}}\]
is the generating function for the overpartitions $\pi$ such that

{\noindent (1)} all parts $\equiv a\pmod A$ of size not exceeding $kA+a$ in $\pi$ are non-overlined.

So, we see that
\[\frac{(-q;q)_{\infty }}{(q;q)_{\infty}}\frac{q^{(r-1)[A{k\choose 2}+ka]}}{(-q^{a};q^{A})_{k+1}}\quad\left(\text{resp. }\frac{(-q;q)_{\infty }}{(q;q)_{\infty}}\frac{q^{(r-1)[A{{k+1}\choose 2}+(k+1)a]}}{(-q^{a};q^{A})_{k+1}}\right)\]
is the generating function for the overpartitions $\pi$ satisfying $\widetilde{mes}_{r,A,a}(\pi)\geq kA+a$ (resp. $\widetilde{mes}_{r,A,a}(\pi)\geq (k+1)A+a$) and the condition (1) above. It yields that
\[\frac{(-q;q)_{\infty }}{(q;q)_{\infty}}\frac{q^{(r-1)[A{k\choose 2}+ka]}}{(-q^{a};q^{A})_{k+1}}-\frac{(-q;q)_{\infty }}{(q;q)_{\infty}}\frac{q^{(r-1)[A{{k+1}\choose 2}+(k+1)a]}}{(-q^{a};q^{A})_{k+1}}\]
is the generating function for the overpartitions $\pi$ satisfying $\widetilde{mes}_{r,A,a}(\pi)=kA+a$ and the condition (1) above. This completes the proof.  \qed

Finally, we give two proofs of Theorem \ref{gen-tildeOmes} with a similar argument in the proofs of Theorem \ref{gen-tildeNmes}.

{\noindent \bf The first proof of Theorem \ref{gen-tildeOmes}.} By definition, we find that for $k\geq 0$, $\pi$ is an overpartition in ${O\widetilde{mes}}_{r,A,a}$ with $\widetilde{mes}_{r,A,a}(\pi)=kA+a$ if and only if
\[f_{\overline{iA+a}}(\pi)=1\text{ and }f_{iA+a}(\pi)=0,\text{ for } 0\leq i\leq k-1,\]
and
\[f_{\overline{kA+a}}(\pi)=0\text{ and }f_{kA+a}(\pi)\leq r-2.\]
It yields that
\begin{align*}
		&\sum_{\pi \in O\widetilde{mes}_{r,A,a} }z^{\widetilde{mes}_{r,A,a}(\pi)}q^{|\pi|}\\
		&=\sum_{k=0}^{\infty } z^{kA+a}q^{[a+(A+a)+\cdots+(k-1)A+a]}(1+q^{kA+a}+\cdots+q^{(r-2)(kA+a)})
		\frac{(q^{a};q^{A})_{k+1}}{(q;q)_{\infty }} \frac{(-q;q)_{\infty }}{(-q^{a};q^{A})_{k+1}} \\
		&=\frac{(-q;q)_{\infty }}{(q;q)_{\infty}} \sum_{k=0}^{\infty }z^{kA+a}q^{A{k\choose	2}+ka}(1-q^{(r-1)(kA+a)})\frac{(q^{a};q^{A})_{k}}{(-q^{a};q^{A})_{k+1}}
		\\
		&=\frac{(-q;q)_{\infty }}{(q;q)_{\infty}}\sum_{k=0}^{\infty }z^{kA+a}\bigg[q^{A{k\choose 2}+ka}\frac{(q^{a};q^{A})_{k}}{(-q^{a};q^{A})_{k+1}}-q^{A{k\choose 2}+ka+(r-1)(kA+a)}\frac{(q^{a};q^{A})_{k}}{(-q^{a};q^{A})_{k+1 }}\bigg]. \nonumber
\end{align*}
The proof is complete.  \qed

{\noindent \bf The second proof of Theorem \ref{gen-tildeOmes}.} For $k\geq 0$, it is clear that
\[\frac{(-q;q)_{\infty}}{(q;q)_{\infty}}q^{A{k\choose 2}+ka}\frac{(q^{a};q^{A})_{k}}{(-q^{a};q^{A})_{k+1}}\]
is the generating function for the overpartitions $\pi$ such that
\begin{itemize}
\item[(1)]  $\widetilde{mes}_{r,A,a}\geq kA+a$;
\item[(2)] there is no overlined part of size $kA+a$ in $\pi$;
\item[(3)] all parts $\equiv a\pmod A$ of size less than $kA+a$ are overlined.
\end{itemize}
Moreover,
\[\frac{(-q;q)_{\infty }}{(q;q)_{\infty}}q^{A{k\choose 2}+ka+(r-1)(kA+a)}\frac{(q^{a};q^{A})_{k}}{(-q^{a};q^{A})_{k+1}}\]
is the generating function for the overpartitions $\pi$ satisfying $\widetilde{mes}_{r,A,a}\geq (k+1)A+a$ and the conditions (2) and (3) above. It yields that
\[\frac{(-q;q)_{\infty}}{(q;q)_{\infty}}q^{A{k\choose 2}+ka}\frac{(q^{a};q^{A})_{k}}{(-q^{a};q^{A})_{k+1}}-\frac{(-q;q)_{\infty }}{(q;q)_{\infty}}q^{A{k\choose 2}+ka+(r-1)(kA+a)}\frac{(q^{a};q^{A})_{k}}{(-q^{a};q^{A})_{k+1}}\]
is the generating function for the overpartitions $\pi$ satisfying $\widetilde{mes}_{r,A,a}=kA+a$ and the conditions (2) and (3) above. This completes the proof.  \qed

\end{document}